\overfullrule=0pt
\centerline {\bf A note on multiple solutions for Kirchhoff-type equations with a Neumann condition}\par
\bigskip
\bigskip
\centerline {BIAGIO RICCERI}\par
\bigskip
\bigskip
\centerline {\it In memory of Professor Franco Giannessi}\par
\bigskip
\bigskip
{\bf Abstract:} Using as a main tool our recent result on the strict minimax inequality proved in [5], in this note we establish a multiplicity theorem
for a problem of the type
$$\cases{-K\left(\int_{\Omega}|\nabla u(x)|^2dx\right)\Delta u = h(x,u) & in $\Omega$\cr & \cr
{{\partial u}\over {\partial\nu}}=0 & on $\partial\Omega$.\cr}$$
\bigskip
\bigskip
{\bf Keywords:} Kirchhoff function; Neumann condition; strict minimax inequality; multiplicity.\par
\bigskip
\bigskip
{\bf 2020 Mathematics Subject Classifications:} 35J20; 35J25.
\bigskip
\bigskip
\bigskip
\bigskip
Let $\Omega\subset {\bf R}^n$ be a smooth bounded domain ($n\geq 2$). We consider $H^1(\Omega)$ endowed with the usual
norm
$$\|u\|=\left(\int_{\Omega}|\nabla u(x)|^2dx+\int_{\Omega}|u(x)|^2dx\right)^{1\over 2}.$$
We denote by ${\cal A}$ the class of all
Carath\'eodory functions $h:\Omega\times {\bf R}\to {\bf R}$ such that
$$\sup_{(x,t)\in \Omega\times {\bf R}}{{|h(x,t)|}\over {|t|^q+1}}<+\infty$$
where $0<q<+\infty$, with $q<{{n+2}\over {n-2}}$ if $n\geq 3$. \par
\smallskip
Given a continuous function $K:[0,+\infty[\to {\bf R}$ and $h\in {\cal A}$ we are interested in the nonlocal Neumann problem
$$\cases{-K\left(\int_{\Omega}|\nabla u(x)|^2dx\right)\Delta u = h(x,u) & in $\Omega$\cr & \cr
{{\partial u}\over {\partial\nu}}=0 & on $\partial\Omega$\cr}$$
where ${{\partial u}\over {\partial\nu}}$ is the outer unit normal to $\partial\Omega$. A weak solution of this propblem is any $u\in H^1(\Omega)$
such that
$$K\left(\int_{\Omega}|\nabla u(x)|^2dx\right)\int_{\Omega}\nabla u(x)\nabla v(x)dx=\int_{\Omega}h(x,u(x))v(x)dx$$
for all $v\in H^1(\Omega)$.  So, by classical results, the weak solutions of the problem agree with
the critical points in $H^1(\Omega)$ of the functional $$u\to {{1}\over {2}}\tilde K\left(\int_{\Omega}|\nabla u(x)|^2dx\right)-\int_{\Omega}H(x,u(x))dx,$$
where $\tilde K(t)=\int_0^tK(s)ds$ and $H(x,t)=\int_0^{t}h(x,s)ds$. Moreover, the functional $u\to \int_{\Omega}H(x,u(x))dx$ turns out to be
sequentially weakly continuous in $H^1(\Omega)$.\par
\smallskip
The aim of this very short note is merely to establish the following result:\par
\medskip
THEOREM 1. - {\it Let $K:[0,+\infty[\to [0,+\infty[$, $f, g:{\bf R}\to {\bf R}$ be three continuous functions satisfying the following
conditions:\par
\noindent
$(a_1)$\hskip 5pt $\int_0^tK(s)ds>0$ for all $t>0$ and $\liminf_{t\to +\infty}K(t)>0$;\par
\noindent
$(a_2)$\hskip 5pt $f\in {\cal A}$, $f$ is odd, $\lim_{|t|\to +\infty}\int_0^tf(s)ds=+\infty$ and $\inf_{t\in {\bf R}}\int_0^tf(s)ds<0$;\par
\noindent
$(a_3)$\hskip 5pt $g\in {\cal A}$, $g$ is even, $\sup_{t\in {\bf R}}\left|\int_0^tg(s)ds\right |<+\infty$ and $\int_0^{t}g(s)ds\neq 0$ for all $ t\in
f^{-1}(0)\cap ]0,+\infty[$.\par
Then, for each $\alpha, \beta\in L^{\infty}(\Omega)$, with $\sup_{\Omega}\alpha<0$, $\int_{\Omega}\beta(x)dx\neq 0$, and for each $r>0$, there exists $\epsilon^*>0$
with the following property: for every $h\in {\cal A}$ such that 
$$\int_{\Omega}\sup_{t\in\bf R}\left(\int_0^th(x,s)ds\right)dx-\int_{\Omega}\inf_{t\in\bf R}\left(\int_0^th(x,s)ds\right)dx<\epsilon^*,$$
 there exists an open interval $A\subset [-r,r]$ such that, for each $\lambda\in A$, the problem
$$\cases{-K\left(\int_{\Omega}|\nabla u(x)|^2dx\right)\Delta u= \alpha(x)f(u)+\lambda\beta(x)g(u)+h(x,u) & in $\Omega$\cr & \cr
{{\partial u}\over {\partial\nu}}=0 & on $\partial\Omega$\cr}\eqno{(P)}$$
has at least two weak solutions.}\par
\medskip
Our proof of Theorem 1 is based on a number of our previous results. For the reader's convenience, we now recall them, not in full generality, but
in particolar forms which are enough for our purposes.\par
\medskip
PROPOSITION A ([3], Proposition 1.2). - {\it Let $\eta:\Omega\times {\bf R}\times {\bf R}^n\to {\bf R}$ be a Carath\'eodory function such that,
for some $c, d>0$, one has
$$c|\theta|^2-d\leq \eta(x,t,\theta)$$
for all $(x,t,\theta)\in \Omega\times {\bf R}\times {\bf R}^n$ and
$$\lim_{|t|\to +\infty}\inf_{(x,\theta)\in \Omega\times {\bf R}^n}\eta(x,t,\theta)=+\infty.$$
Then, in $H^1(\Omega)$, one has
$$\lim_{\|u\|\to +\infty}\int_{\Omega}\eta(x,u(x),\nabla u(x))=+\infty.$$}
\medskip
THEOREM A ([5], Theorem 2.2). - {\it Let $X$ be a reflexive real Banach space, $Q, \psi: X\to {\bf R}$, $r>0$. Assume that $Q$ is even, that
$\psi$ is odd and that, for each $\lambda\in [-r,r]$, the functional $Q-\lambda\psi$ is sequentially weakly lower semicontinuous and
coercive. Finally, suppose that there is no global minimum of $Q$ at which $\psi$ vanishes.\par
Then, one has
$$\sup_{|\lambda|\leq r}\inf_{x\in X}(Q(x)-\lambda\psi(x))<\inf_{x\in X}\sup_{|\lambda|\leq r}(Q(x)-\lambda\psi(x)).$$}
\medskip
PROPOSITION B ([4], Theorem 1). - {\it Let $X, Y$ be two non-empty sets and let $k:X\times Y\to {\bf R}$ be such that
$$\sup_Y\inf_Xk<\inf_X\sup_Yk.$$
Then, for each $\varphi:X\to {\bf R}$ satisfying
$$\sup_X\varphi-\inf_X\varphi<\inf_X\sup_Yk-\sup_Y\inf_Xk,$$
one has
$$\sup_Y\inf_X(k+\varphi)<\inf_X\sup_Y(k+\varphi).$$}
\medskip
THEOREM B ([2], Theorem B). - {\it . Let $X$ be a reflexive real Banach space, $I\subseteq {\bf R}$ an interval and $k:X\times I\to {\bf R}$
a function such that $k(x,\cdot)$ is concave
in $I$ for all $x\in X$, while $k(\cdot,\lambda)$ is continuous, coercive and sequentially weakly lower semicontinuous in $X$ for all $\lambda\in I$. Further,
assume that
$$\sup_Y\inf_Xk<\inf_X\sup_Yk.$$
Then, there exists an open interval $A\subseteq I$, such that, for each $\lambda\in A$, the function $k(\cdot,\lambda)$ has at least two local minima.}\par
\medskip
Now, we can give the\par
\medskip
{\it Proof of Theorem 1}. Consider the functionals $Q, \psi:H^1(\Omega)\to {\bf R}$ defined by
$$Q(u)={{1}\over {2}}\tilde K\left(\int_{\Omega}|\nabla u(x)|^2dx\right)-\int_{\Omega}\alpha(x)F(u(x))dx$$
$$\psi(u)=-\int_{\Omega}\beta(x)G(u(x))dx,$$
where $\tilde K(t)=\int_0^tK(s)ds$,  $F(t)=\int_0^tf(s)ds$, $G(t)=\int_0^sg(s)ds$. In view of $(a_2)$ and since $\sup_{\Omega}\alpha<0$, we
have
$$\lim_{|t|\to +\infty}\inf_{x\in \Omega}-\alpha(x)F(t)=+\infty.\eqno{(1)}$$
On the other hand, in view of $(a_1)$, there exist $c, d>0$ so that
$$ct-d\leq {{1}\over {2}}\tilde K(t)-\alpha(x)F(s)\eqno{(2)}$$
for all $t\geq 0, s\in {\bf R}, x\in\Omega$. At this point, $(1)$ and $(2)$ allow us to apply Proposition A, obtaining
$$\lim_{\|u\|\to +\infty}Q(u)=+\infty.\eqno{(3)}$$
Of course, we have
$$\inf_{H^1(\Omega)}Q= -\int_{\Omega}\alpha(x)dx\inf_{{\bf R}}F.$$
Now, let $u\in H^1(\Omega)$ satisfy
$$Q(u)= -\int_{\Omega}\alpha(x)dx\inf_{{\bf R}}F.$$
That is
$${{1}\over {2}}\tilde K\left(\int_{\Omega}|\nabla u(x)|^2dx\right)=\int_{\Omega}\alpha(x)\left(F(u(x))-\inf_{{\bf R}}F\right)dx.\eqno{(4)}$$
Since the left-hand side of $(4)$ is non-negative, while the right-hand side is non-positive, we infer that
$\int_{\Omega}|\nabla u(x)|^2dx=0$ (by $(a_1)$) and $F(u(x))=\inf_{{\bf R}}F$ a.e. in $\Omega$. Consequently, since $\Omega$ is connected,
$u$ is constant. Viceversa, if $\tilde t\in {\bf R}$ satisfies $F(\tilde t)=\inf_{{\bf R}}F$, then the constant function $u(x)=\tilde t$ satisfies $(4)$.
That is to say, the global minima of the functional $Q$ are exactly the constant functions whose values are the global minima of $F$ in ${\bf R}$.
Let $\tilde t\in {\bf R}$ be a global minimum of $F$. Then, we have $f(\tilde t)=0$ and, by $(a_2)$, $\tilde t\neq 0$. Consequently, by $(a_3)$ and
since $\int_{\Omega}\beta(x)dx\neq 0$, there is no global minimum of $Q$ at which $\psi$ vanishes. Of course, $Q$ is even and $\psi$ is odd.
Further, notice that, since $\tilde K$ is continuous and non-decreasing, the functional $u\to \tilde K\left (\int_{\Omega}|\nabla u(x)|^2dx\right)$
is weakly lower semicontinuous in $H^1(\Omega)$. Consequently, for each $\lambda\in [-r,r]$, the functional $Q-\lambda\psi$ is
sequentially weakly lower semicontinuous and coercive (by $(3)$ and $(a_3)$). Therefore, we can apply Theorem A, obtaining
$$\sup_{|\lambda|\leq r}\inf_{u\in H^1(\Omega)}(Q(u)-\lambda\psi(u))<\inf_{u\in H^1(\Omega)}\sup_{|\lambda|\leq r}(Q(u)-\lambda\psi(u)).$$
Now, put
$$\epsilon^*=\inf_{u\in H^1(\Omega)}\sup_{|\lambda|\leq r}(Q(u)-\lambda\psi(u))-\sup_{|\lambda|\leq r}\inf_{u\in H^1(\Omega)}(Q(u)-\lambda\psi(u)).$$
Fix any $h\in {\cal A}$ such that
$$\int_{\Omega}\sup_{t\in\bf R}\left(\int_0^th(x,s)ds\right)dx-\int_{\Omega}\inf_{t\in\bf R}\left(\int_0^th(x,s)ds\right)dx<\epsilon^*.$$
Set $H(x,t)=\int_0^th(x,s)ds$. Consider the functional $\varphi:H^1(\Omega)\to {\bf R}$ defined by
$$\varphi(u)=-\int_{\Omega}H(x,u(x))dx.$$
Since
$$\int_{\Omega}\inf_{t\in {\bf R}}H(x,t)dx\leq \varphi(u)\leq \int_{\Omega}\sup_{t\in {\bf R}}H(x,t)dx,$$
we have
$$\sup_{H^1(\Omega)}\varphi-\inf_{H^1(\Omega)}\varphi\leq \int_{\Omega}\sup_{t\in {\bf R}}H(x,t)dx-\int_{\Omega}\inf_{t\in {\bf R}}H(x,t)dx<
\epsilon^*.$$
At this point, we can apply Proposition B, obtaining
$$\sup_{|\lambda|\leq r}\inf_{u\in H^1(\Omega)}(Q(u)-\lambda\psi(u)-\varphi(u))<\inf_{u\in H^1(\Omega)}\sup_{|\lambda|\leq r}(Q(u)-\lambda\psi(u))-\varphi(u)).$$
Finally, applying Theorem B, we obtain the existence of an open interval $A\subseteq [-r,r]$ such that, for each $\lambda\in A$, the functional
$Q-\lambda\psi-\varphi$ has at least two local minima in $H^1(\Omega)$: they are the claimed weak solutions of problem $(P)$.\hfill $\bigtriangleup$
\medskip
REMARK 1. - We have not found in the literature any result close enough to Theorem 1 in order to do a proper comparison. We refer to [1] for a
very recent contribution to the topic.\par
\medskip
This note leaves the following hard problem open:\par
\medskip
PROBLEM 1. - Assume that the assumptions of Theorem 1 are satisfied. Prove or disprove the following: \par
For each $\alpha, \beta\in L^{\infty}(\Omega)$, with $\sup_{\Omega}\alpha<0$, $\int_{\Omega}\beta(x)dx\neq 0$, there exist $\epsilon^*, \rho^*>0$
with the following property: for every $h\in {\cal A}$ such that 
$$\int_{\Omega}\sup_{t\in\bf R}\left(\int_0^th(x,s)ds\right)dx-\int_{\Omega}\inf_{t\in\bf R}\left(\int_0^th(x,s)ds\right)dx<\epsilon^*,$$
and for every $\lambda\in]0,\rho^*[$, the problem
$$\cases{-K\left(\int_{\Omega}|\nabla u(x)|^2dx\right)\Delta u= \alpha(x)f(u)+\lambda\beta(x)g(u)+h(x,u) & in $\Omega$\cr & \cr
{{\partial u}\over {\partial\nu}}=0 & on $\partial\Omega$\cr}$$
has at least two weak solutions.\par
\bigskip
\bigskip
{\bf Acknowledgements.} This work has been funded by the European Union - NextGenerationEU Mission 4 - Component 2 - Investment 1.1 under the Italian Ministry of University and Research (MUR) programme "PRIN 2022" - grant number 2022BCFHN2 - Advanced theoretical aspects in PDEs and their applications - CUP: E53D23005650006. The author has also been supported by the Gruppo Nazionale per l'Analisi Matematica, la Probabilit\`a e 
le loro Applicazioni (GNAMPA) of the Istituto Nazionale di Alta Matematica (INdAM) and by the Universit\`a degli Studi di Catania, PIACERI 2024-2026, Linea di intervento 2, Progetto ”PAFA”.
\vfill\eject
\centerline {\bf References}\par
\bigskip
\bigskip
\noindent
[1]\hskip 5pt F. BORER, M. T. O. PIMENTA and P. WINKERT, 
{\it Degenerate Kirchhoff problems with nonlinear Neumann boundary condition}, J. Funct. Anal. {\bf 289} (2025), no. 4, Paper No. 110933.\par
\smallskip
\noindent
[2]\hskip 5pt B. RICCERI, {\it A three critical points theorem revisited}, Nonlinear Anal., {\bf 70} (2009), 3084-3089.\par
\smallskip
\noindent
[3]\hskip 5pt B. RICCERI, {\it Miscellaneous applications of certain minimax theorems II}, Acta Math. Vietnam., {\bf 45} (2020), 515-524.\par
\smallskip
\noindent
[4]\hskip 5pt B. RICCERI, {\it Addendum to  ``A more complete version of a minimax theorem"}, Appl. Anal. Optim., {\bf 6} (2022), 195-197.\par
\smallskip
\noindent
[5]\hskip 5pt B. RICCERI, {\it Multiple critical points in closed sets via minimax theorems}, Optimization, 1-15.\par
 \noindent
https://doi.org/10.1080/02331934.2025.2457549.
\bigskip
\bigskip
\bigskip
\bigskip
Department of Mathematics and Informatics\par
University of Catania\par
Viale A. Doria 6\par
95125 Catania, Italy\par
{\it e-mail address}: ricceri@dmi.unict.it

\bye